\input amstex
%
%
%
%
\documentstyle{amsppt}
\leftheadtext{}
\rightheadtext{}
\magnification=1000
\baselineskip=25pt
\hoffset=0.12truein
\voffset=0.5truein
\NoBlackBoxes
\pagewidth{15.5truecm}
\pageheight{20truecm}
\topmatter
\title $k$-Hyponormality of Finite Rank Perturbations of
Unilateral Weighted Shifts\endtitle
\author Ra\'ul E. Curto and Woo Young Lee\endauthor
\address Department of Mathematics, University of Iowa,
Iowa City, IA 52242\endaddress
\email rcurto\@math.uiowa.edu\endemail
\address Department of Mathematics, Seoul National University, Seoul 151-742, Korea
\endaddress
\email wylee\@math.snu.ac.kr\endemail
\thanks The work of the first-named author was partially supported
by NSF research grants DMS-9800931 and DMS-0099357.
\endthanks
\thanks The work of the second-named author was partially
supported by a grant (R14-2003-006-01001-0) from the Korea Science
and Engineering Foundation.
\endthanks
\abstract
In this paper
we explore finite rank perturbations of unilateral weighted shifts $W_\alpha$.
First, we prove that the
subnormality of $W_\alpha$  is
never stable
under nonzero finite rank pertrubations unless
the perturbation
occurs at the zeroth weight.
Second,
we establish that 2-hyponormality implies positive quadratic
hyponormality, in the sense that the Maclaurin coefficients of
$D_n(s):=\text{det}\,P_n\,[(W_\alpha+sW_\alpha^2)^*,\,
W_\alpha+s W_\alpha^2]\,P_n$
are nonnegative, for every $n\ge 0$, where $P_n$ denotes the orthogonal
projection onto the basis vectors $\{e_0,\cdots,e_n\}$.
Finally, for $\alpha$ strictly increasing and $W_\alpha$ 2-hyponormal,
we show that
for a small finite-rank perturbation $\alpha^\prime$ of $\alpha$, the shift
$W_{\alpha^\prime}$ remains quadratically hyponormal.
\endabstract
\endtopmatter
\document
\head
1. Introduction
\endhead
\bigskip

Let $\Cal{H}$ and $\Cal{K}$ be complex Hilbert spaces,
let $\Cal{L(H,K)}$ be the set of
bounded linear operators from $\Cal{H}$ to $\Cal{K}$ and write
$\Cal{L(H)}:=\Cal{L(H,H)}$.
An operator $T\in\Cal{L(H)}$ is said to be normal if $T^*T=TT^*$,
hyponormal if $T^*T\ge TT^*$, and subnormal if $T=N\vert_{\Cal{H}}$,
where $N$ is normal on some Hilbert space $\Cal{K}\supseteq \Cal{H}$.
If $T$ is subnormal then $T$ is also hyponormal. Recall that given a
bounded sequence of positive numbers $\alpha:\alpha_0,\alpha_1,\cdots$
(called {\it weights}), the {\it (unilateral) weighted shift} $W_\alpha$
associated with $\alpha$ is the operator on $\ell^2(\Bbb{Z}_+)$ defined by
$W_\alpha e_n:=\alpha_n e_{n+1}$ for all $n\ge 0$, where
$\{e_n\}_{n=0}^\infty$ is the canonical orthonormal basis for $\ell^2$.
It is straightforward to check that $W_\alpha$ can never be {\it normal},
and that $W_\alpha$ is {\it hyponormal} if and only if
$\alpha_n\le \alpha_{n+1}$ for all $n\ge 0$.
The Bram-Halmos criterion for subnormality states that an operator
$T$ is subnormal if
and only if
$$
\sum_{i,j} (T^i x_j, T^j x_i)\ge 0
$$
for all finite collections $x_0,x_1,\cdots,x_k\in\Cal{H}$
([2],[4, II.1.9]).
It is easy to see that this is equivalent to the following positivity test:
$$
\pmatrix
I&T^*&\hdots& T^{*k}\\
T& T^*T& \hdots& T^{*k}T\\
\vdots&\vdots& \ddots& \vdots\\
T^k& T^*T^k& \hdots& T ^{*k}T^k\endpmatrix\ge 0\qquad\text{(all $k\ge 1$)}.
\tag 1.1
$$
Condition (1.1) provides a measure of the gap
between hyponormality and subnormality. In fact, the positivity condition
(1.1) for $k=1$ is equivalent to the hyponormality of $T$, while
subnormality requires the validity of (1.1) for all $k$.
Let $[A,B]:=AB-BA$ denote the commutator of two operators $A$ and $B$,
and define $T$ to be {\it $k$-hyponormal} whenever the $k\times k$
operator matrix
$$
M_k(T):=([T^{*j},T^i])_{i,j=1}^k
\tag 1.2
$$
is positive.
An application of the Choleski algorithm for operator matrices
shows that the positivity of (1.2) is equivalent to
the positivity of
the $(k+1)\times (k+1)$ operator matrix in
(1.1); the Bram-Halmos criterion can be then rephrased as
saying that $T$ is subnormal if and only if $T$ is $k$-hyponormal for
every $k\ge 1$ ([16]).

Recall ([1],[16],[5])
that $T\in\Cal{L(H)}$ is said to be
{\it weakly $k$-hyponormal} if
$$
LS(T,T^2,\cdots,T^k):=\left\{\sum_{j=1}^k \alpha_jT^j: \alpha=
(\alpha_1,\cdots,\alpha_k)\in {\Bbb C}^k\right\}
$$
consists entirely of hyponormal operators,
or equivalently, $M_k(T)$ is {\it weakly positive}, i.e., ([16])
$$
(M_k(T)\pmatrix \lambda_0 x\\ \vdots\\ \lambda_{k}x\endpmatrix,\ \pmatrix
\lambda_0 x\\ \vdots\\ \lambda_{k}x\endpmatrix)\ \ge 0\qquad
\text{for $x\in \Cal{H}$ and $\lambda_0,\cdots,\lambda_k \in{\Bbb C}$}.
\tag 1.3
$$
If $k=2$ then $T$ is said to be {\it quadratically hyponormal},
and if $k=3$ then $T$ is said to be {\it cubically hyponormal}.
Similarly, $T\in\Cal{L(H)}$ is said to be {\it polynomially hyponormal}
if $p(T)$ is hyponormal for every polynomial $p\in {\Bbb C}[z]$.
It is known that $k$-hyponormal $\Rightarrow$ weakly $k$-hyponormal,
but the converse is not true in general.

The classes of (weakly) $k$-hyponormal operators have been studied
in an attempt to bridge the gap between subnormality and hyponormality
([7],[8],[10],[11],[12],[14],[16],[19],[22]). The study of this gap
has been only partially successful.
For example, such a gap is not yet well described for Toeplitz operators
on the Hardy space of the unit circle; in fact, even subnormality for
Toeplitz operators has not been characterized (cf.[20],[6]).
For weighted shifts, positive results appear in
[7] and [12], although no concrete example of a
weighted shift which is polynomially hyponormal and not subnormal
has yet been found (the existence of such weighted
shifts was established in [17] and [18]).
\medskip

In the present paper we renew our efforts to help describe the above
mentioned gap between subnormality and hyponormality, with particular
emphasis on polynomial hyponormality. We
focus on the class of unilateral weighted shifts, and initiate a study of
how the above mentioned notions behave
under finite perturbations of the weight sequence.
We first obtain three concrete results:

(i) the subnormality of $W_\alpha$ is never stable under nonzero
finite rank perturbations unless the perturbation is confined to the
zeroth weight (Theorem 2.1);

(ii) 2-hyponormality implies {\it positive quadratic
hyponormality}, in the sense that the Maclaurin coefficients of
$D_n(s):=\text{det}\,P_n\,[(W_\alpha+sW_\alpha^2)^*,\,
W_\alpha+s W_\alpha^2]\,P_n$
are nonnegative, for every $n\ge 0$, where $P_n$ denotes the orthogonal
projection onto the basis vectors $\{e_0,\cdots,e_n\}$ (Theorem 2.2); and

(iii) if $\alpha$ is strictly increasing and
$W_\alpha$ is $2$-hyponormal then for $\alpha^\prime$ a small
perturbation of $\alpha$, the shift $W_{\alpha^\prime}$ remains
positively quadratically hyponormal (Theorem 2.3).

Along the way we establish two related results, each of independent interest:

(iv) an integrality criterion for a subnormal weighted shift to
have an $n$-step
subnormal extension (Theorem 6.1); and

(v) a proof that the sets of $k$-hyponormal and weakly $k$-hyponormal
operators are closed in the strong operator topology (Proposition 6.7).
\bigskip

\head
2. Statement of Main Results
\endhead
\bigskip

C. Berger's characterization of
subnormality for unilateral weighted shifts (cf. [21],[4, III.8.16])
states that $W_\alpha$ is subnormal
if and only if there exists a Borel probability measure $\mu$ (the so-called Berger measure of $W_\alpha$)
supported in $[0,||W_\alpha||^2]$, with $||W_\alpha||^2\in\text{supp}\,\mu$,
such that
$$
\gamma_n=\int t^n d\mu(t)\quad\text{for all}\ n\ge 0.
$$
Given an initial segment of weights $\alpha:\alpha_0,\cdots\alpha_m$,
the sequence $\hat\alpha\in \ell^\infty({\Bbb Z}_+)$ such that
$\hat\alpha_i=\alpha_i\ (i=0,\cdots,m)$ is said to be
{\it recursively generated} by $\alpha$ if there exist $r\ge 1$
and $\varphi_0,\cdots,\varphi_{r-1}\in {\Bbb R}$ such that
$$
\gamma_{n+r}=\varphi_0\gamma_n+\cdots+\varphi_{r-1}\gamma_{n+r-1}
\quad\text{(all $n\ge 0$)},
\tag 2.1
$$
where $\gamma_0:=1$, $\gamma_n:=\alpha_0^2\cdots\alpha_{n-1}^2$ ($n\ge 1$).
In this case $W_{\hat\alpha}$ with weights $\hat\alpha$ is said to
be {\it recursively generated}. If we let
$$
g(t):=t^r-\left(\varphi_{r-1}t^{r-1}+\cdots+\varphi_0\right),
\tag 2.2
$$
then $g$ has $r$ distinct real roots $0\le s_0<\cdots<s_{r-1}$
([11, Theorem 3.9]). Let
$$
V:=\pmatrix 1&1&\hdots&1\\
s_0&s_1&\hdots&s_{r-1}\\
\vdots&\vdots&&\vdots\\
s_0^{r-1}&s_1^{r-1}&\hdots&s_{r-1}^{r-1}\endpmatrix
$$
and let
$$
\pmatrix
\rho_0\\ \vdots\\ \rho_{r-1}\endpmatrix:=V^{-1}\pmatrix \gamma_0\\
\vdots\\ \gamma_{r-1}\endpmatrix.
$$
If the associated recursively generated weighted shift $W_{\hat\alpha}$
is subnormal then its Berger measure is
of the form
$$
\mu:=\rho_0\delta_{s_0}+\cdots+\rho_{r-1}\delta_{r-1}.
$$
For example, given $\alpha_0<\alpha_1<\alpha_2$,
$W_{(\alpha_0,\alpha_1,\alpha_2)^\wedge}$
is the recursive
weighted shift whose weights are calculated according to the recursive
relation
$$
\alpha_{n+1}^2=\varphi_1+\varphi_0\frac{1}{\alpha_n^2},
\tag 2.3
$$
where
$$
\varphi_0=-\frac{\alpha_0^2\alpha_1^2(\alpha_2^2-\alpha_1^2)}
{\alpha_1^2-\alpha_0^2}\quad\text{and}\quad
\varphi_1=\frac{\alpha_1^2(\alpha_2^2-\alpha_0^2)}{\alpha_1^2-\alpha_0^2}.
\tag 2.4
$$
In this case, $W_{(\alpha_0,\alpha_1,\alpha_2)^\wedge}$ is
subnormal with $2$--atomic Berger measure.
Let  $W_{x\,(\alpha_0,\alpha_1,\alpha_2)^\wedge}$ denote
the weighted shift whose weight sequence consists of the initial weight
$x$ followed by the weight sequence of $W_{(\alpha_0,\alpha_1,
\alpha_2)^\wedge}$.

By the Density Theorem ([11, Theorem 4.2 and Corollary 4.3]),
we know that if $W_\alpha$ is a subnormal weighted shift with weights
$\alpha=\{\alpha_n\}$ and $\epsilon>0$, then there exists a nonzero
compact operator $K$ with $||K||<\epsilon$ such that $W_\alpha+K$ is
a recursively generated subnormal weighted shift;
in fact $W_\alpha+K=W_{\widehat{\alpha^{(m)}}}$ for some $m\ge 1$, where
$\alpha^{(m)}:\alpha_0,\cdots,\alpha_m$.
The following result shows that $K$ cannot generally be taken
to be finite rank.
\smallskip

\proclaim{Theorem 2.1 (Finite Rank Perturbations of Subnormal Shifts)}
If $W_\alpha$ is a subnormal weighted shift then there exists no
nonzero finite rank operator $F (\ne cP_{\{e_
0\}})$ such that
$W_\alpha+F$ is a subnormal weighted shift.
Concretely, suppose $W_\alpha$ is a subnormal weighted shift
with weight sequence $\alpha=\{\alpha_n\}_{n=0}^\infty$ and assume
$\alpha^\prime=\{\alpha_n^\prime\}$ is a nonzero perturbation of
$\alpha$ in a finite number of weights except the initial
weight; then $W_{\alpha^\prime}$ is not subnormal.
\endproclaim
\bigskip

We next consider the self-commutator
$[(W_\alpha+s\,W_\alpha^2)^*, W_\alpha+s\,W_\alpha^2]$.
Let $W_\alpha$ be a hyponormal weighted shift. For $s\in \Bbb{C}$,
we write
$$
D(s):= [(W_\alpha+s\,W_\alpha^2)^*, W_\alpha+s\,W_\alpha^2]
$$
and we let
$$
D_n(s):=P_n[(W_\alpha+s\,W_\alpha^2)^*, W_\alpha+s\,W_\alpha^2]P_n
=\pmatrix
q_0&\bar r_0&0&\hdots&0&0\\
r_0&q_1&\bar r_1&\hdots&0&0\\
0&r_1&q_2&\hdots&0&0\\
\vdots&\vdots&\vdots&\ddots&\vdots&\vdots\\
0&0&0&\hdots&q_{n-1}&\bar r_{n-1}\\
0&0&0&\hdots&r_{n-1}&q_n\endpmatrix,
\tag 2.5
$$
where $P_n$ is the orthogonal projection onto the subspace
generated by $\{e_0,\cdots,e_n\}$,
$$
\cases
q_n:=u_n+|s|^2v_n\\
r_n:=s\sqrt {w_n}\\
u_n:=\alpha_n^2-\alpha_{n-1}^2\\
v_n:=\alpha_n^2\alpha_{n+1}^2-\alpha_{n-1}^2\alpha_{n-2}^2\\
w_n:=\alpha_n^2(\alpha_{n+1}^2-\alpha_{n-1}^2)^2,
\endcases
\tag 2.6
$$
and, for notational convenience, $\alpha_{-2}=\alpha_{-1}=0$.
Clearly, $W_\alpha$ is quadratically hyponormal if and only
if $D_n(s)\ge 0$ for all $s\in\Bbb{C}$ and all $n\ge 0$.
Let $d_n(\cdot):=\text{det}\,(D_n(\cdot))$. Then $d_n$
satisfies the following 2--step recursive formula:
$$
d_0=q_0,\quad d_1=q_0q_1-|r_0|^2,\quad d_{n+2}=
q_{n+2}d_{n+1}-|r_{n+1}|^2d_n.
\tag 2.7
$$
If we let $t:=|s|^2$, we observe that $d_n$ is a polynomial
in $t$ of degree $n+1$, and if we write $d_n\equiv \sum_{i=0}^{n+1}c(n,i)t^i$,
then the coefficients $c(n,i)$ satisfy a double-indexed recursive formula,
namely
$$
\split
c(n+2,i)&=u_{n+2}\,c(n+1,i)+v_{n+2}\,c(n+1,i-1)-w_{n+1}\,c(n,i-1),\\
c(n,0)&=u_0\cdots u_n,\quad c(n,n+1)=v_0\cdots v_n,\quad c(1,1)=
u_1v_0+v_1u_0-w_0
\endsplit
\tag 2.8
$$
($n\ge 0$, $i\ge 1$). We say that $W_\alpha$
is {\it positively quadratically hyponormal} if $c(n,i)\ge 0$
for every $n\ge 0,\ 0\le i\le n+1$ (cf. [9]).
Evidently, positively
quadratically hyponormal $\Longrightarrow$ quadratically
hyponormal. The converse, however, is not true in general (cf. [3]).
\medskip

The following theorem establishes a useful relation between $2$-hyponormality
and positive quadratic hyponormality.

\proclaim{Theorem 2.2} Let $\alpha\equiv\{\alpha_n\}_{n=0}^\infty$
be a weight sequence and assume that $W_\alpha$ is $2$-hyponormal.
Then $W_\alpha$ is positively quadratically hyponormal.
More precisely, if $W_\alpha$ is 2-hyponormal then
$$
c(n,i)\ \ge\ v_0\cdots v_{i-1}u_i\cdots u_n\qquad\text{($n\ge 0,\
0\le i\le n+1$)}.
\tag 2.9
$$
In particular, if $\alpha$ is strictly increasing and $W_\alpha$ is $2$-hyponormal then
the Maclaurin coefficients of $d_n(t)$ are positive for all
$n\ge 0$.
\endproclaim
\bigskip

If $W_\alpha$ is a weighted shift with weight sequence
$\alpha=\{\alpha_n\}_{n=0}^\infty$, then the {\it moments} of $W_\alpha$
are usually defined by $\beta_0:=1,\ \beta_{n+1}:=\alpha_n\beta_n$
($n\ge 0$) [23]; however, we prefer to
reserve this term for the sequence
$\gamma_n:=\beta_n^2$ ($n\ge 0$).
A criterion for $k$-hyponormality can be given in terms
of these moments ([7, Theorem 4]): if we build a $(k+1)\times (k+1)$
Hankel matrix $A(n;k)$ by
$$
A(n;k):=\pmatrix
\gamma_n&\gamma_{n+1}&\hdots&\gamma_{n+k}\\
\gamma_{n+1}&\gamma_{n+2}&\hdots&\gamma_{n+k+1}\\
\vdots&\vdots&&\vdots\\
\gamma_{n+k}&\gamma_{n+k+1}&\hdots&\gamma_{n+2k}\endpmatrix\quad
\text{($n\ge0$)},
\tag 2.10
$$
then
$$
\text{$W_\alpha$ is $k$-hyponormal}\ \Longleftrightarrow\ A(n;k)\ge
0\quad\text{($n\ge 0$)}.
\tag 2.11
$$
In particular, for $\alpha$ strictly increasing,
$W_\alpha$ is 2-hyponormal if and only if
$$
\text{det} \pmatrix \gamma_n&\gamma_{n+1}&\gamma_{n+2}\\
\gamma_{n+1}&\gamma_{n+2}&\gamma_{n+3}\\
\gamma_{n+2}&\gamma_{n+3}&\gamma_{n+4}\endpmatrix\ge 0\quad
\text{($n\ge 0$)}.
\tag 2.12
$$
\medskip

One might conjecture that if $W_\alpha$ is a $k$-hyponormal
weighted shift whose weight sequence is strictly increasing then
$W_\alpha$ remains weakly $k$-hyponormal under a small perturbation of
the weight sequence. We will show below that
this is true for $k=2$ (Theorem 2.3).
\bigskip

In [12, Theorem 4.3], it was shown that the gap
between 2-hyponormality and quadratic hyponormality
can be detected by unilateral shifts with a weight sequence
$\alpha:\sqrt{x},(\sqrt{a},\sqrt{b},\sqrt{c})^\wedge$.
In particular, there exists a maximum value $H_2\equiv
H_2(a,b,c)$ of $x$ that makes $W_{\sqrt{x},(\sqrt{a},\sqrt{b},\sqrt{c})^\wedge}$
2-hyponormal; $H_2$ is called the {\it modulus} of 2-hyponormality (cf. [12]).
Any value of $x>H_2$ yields a non-2-hyponormal weighted shift.
However, if $x-H_2$ is small enough,
$W_{\sqrt{x},(\sqrt{a},\sqrt{b},\sqrt{c})^\wedge}$
is still quadratically hyponormal.
The following theorem shows that, more generally,
for finite rank perturbations of
weighted shifts with strictly increasing weight sequences,
there always exists a gap between 2-hyponormality and
quadratic hyponormality.

\proclaim{Theorem 2.3 (Finite Rank Perturbations of 2-hyponormal
Shifts)}
Let $\alpha=\{\alpha_n\}_{n=0}^\infty$ be a strictly increasing
weight sequence.
If $W_\alpha$ is $2$-hyponormal then
$W_\alpha$ remains positively quadratically hyponormal under a small
nonzero finite
rank perturbation of $\alpha$.
\endproclaim
\bigskip


\head 3. Proof of Theorem 2.1
\endhead
\bigskip

\demo{Proof of Theorem 2.1}
It suffices to show that if $T$ is a weighted shift whose restriction
to $\bigvee\{e_n,e_{n+1},\cdots\}$ ($n\ge 2$) is subnormal then there
is at most one $\alpha_{n-1}$ for which $T$ is subnormal.

Let $W:=T\vert_{\bigvee\{e_{n-1},e_n,e_{n+1},\cdots\}}$ and $S:=
T\vert_{\bigvee\{e_n,e_{n+1},\cdots\}}$, where $n\ge 2$. Then $W$ and
$S$ have weights $\alpha_k(W):=\alpha_{k+n-1}$
and $\alpha_k(S):=\alpha_{k+n}$ ($k\ge 0$).
Thus the corresponding moments are related by the equation
$$
\gamma_k(S)=\alpha_n^2\cdots\alpha_{n+k-1}^2=\frac{\gamma_{k+1}(W)}
{\alpha_{n-1}^2}.
$$
We now adapt the proof of [7, Proposition 8].
Suppose $S$ is subnormal with associated Berger measure $\mu$.
Then $\gamma_k(S)=\int_0^{||T||^2}t^{k}\,d\mu$. Thus $W$ is subnormal
if and only if there exists a probability measure $\nu$ on $[0,||T||^2]$
such that
$$
\frac{1}{\alpha_{n-1}^2}\int_{0}^{||T||^2} t^{k+1}\,d\nu(t)=
\int_{0}^{||T||^2} t^{k}\,d\mu(t)\quad\text{for all}\ k\ge 0,
$$
which readily implies that $t\,d\nu=\alpha_{n-1}^2\,d\mu$. Thus
$W$ is subnormal if and only if the formula
$$
d\nu:=\lambda\cdot\delta_0+\frac{\alpha_{n-1}^2}{t} d\mu
\tag 3.1
$$
defines a probability measure for some $\lambda\ge 0$,
where $\delta_0$ is the point mass at the origin.
In particular $\frac{1}{t}\in L^1(\mu)$ and $\mu(\{0\})=0$
whenever $W$ is subnormal. If we repeat the above argument for $W$
and $V:=T\vert_{\bigvee\{e_{n-2},e_{n-1},\cdots\}}$, then we should
have that $\nu(\{0\})=0$ whenever $V$ is subnormal.
Therefore we can conclude that if $V$ is subnormal
then $\lambda=0$, and hence
$$
d\nu=\frac{\alpha_{n-1}^2}{t} d\mu.
\tag 3.2
$$
Thus we have
$$
1=\int_0^{||T||^2} d\nu(t)=\alpha_{n-1}^2\int_0^{||T||^2}\frac{1}{t} d\mu(t),
$$
so that
$$
\alpha_{n-1}^2=\left(\int_0^{||T||^2}\frac{1}{t} d\mu(t)\right)^{-1},
\tag 3.3
$$
which implies that $\alpha_{n-1}$ is determined uniquely by
$\{\alpha_n,\alpha_{n+1},\cdots\}$ whenever $T$ is subnormal.
This completes the proof.
\hfill$\square$
\enddemo
\bigskip

Theorem 2.1 says that a nonzero finite rank perturbation of a
subnormal shift is never subnormal unless the perturbation
occurs at the initial weight. However, this is not the case for
$k$-hyponormality. To see this we use a close relative of
the Bergman shift $B_+$ (whose weights are given by
$\alpha=\{\sqrt{\frac{n+1}{n+2}}\}_{n=0}^\infty$);
it is well known that $B_+$ is subnormal.

\proclaim{Example 3.1} For $x>0$, let $T_x$ be the weighted shift
whose weights are given by
$$
\alpha_0:=\sqrt{\frac{1}{2}},\ \ \alpha_1:=\sqrt{x},\ \ \text{and}\ \
\alpha_n:=\sqrt{\frac{n+1}{n+2}}\ (n\ge 2).
$$
Then we have:
\roster
\item"(i)" $T_x$ is subnormal $\Longleftrightarrow$ $x=\frac{2}{3}$;
\item"(ii)" $T_x$ is 2-hyponormal $\Longleftrightarrow$ $\frac{63-
\sqrt{129}}{80}\le x\le \frac{24}{35}$.
\endroster
\endproclaim

\demo{Proof}
Assertion (i) follows from Theorem 2.1.
For assertion (ii) we use (2.12): $T_x$ is 2-hyponormal
if and only if
$$
\text{det}\pmatrix
1&\frac{1}{2}&\frac{1}{2}x\\
\frac{1}{2}&\frac{1}{2}x&\frac{3}{8}x\\
\frac{1}{2}x&\frac{3}{8}x&\frac{3}{10}x\endpmatrix\ge 0\quad
\text{and}\quad \text{det}
\pmatrix
\frac{1}{2}&\frac{1}{2}x&\frac{3}{8}x\\
\frac{1}{2}x&\frac{3}{8}x&\frac{3}{10}x\\
\frac{3}{8}x&\frac{3}{10}x&\frac{1}{4}x\endpmatrix\ge 0,
$$
or equivalently, $\frac{63-\sqrt{129}}{80}\le x\le \frac{24}{35}$.
\hfill$\square$
\enddemo
\bigskip

For perturbations of recursive subnormal shifts
of the form $W_{(\sqrt{a},\sqrt{b},\sqrt{c})^\wedge}$, subnormality
and 2-hyponormality coincide.

\proclaim{Theorem 3.2} Let $\alpha=\{\alpha_n\}_{n=0}^\infty$
be recursively generated by $\sqrt{a},\sqrt{b},\sqrt{c}$. If $T_x$
is the weighted shift whose weights are given by
$\alpha_x:\alpha_0,\cdots,\alpha_{j-1},\sqrt{x},\alpha_{j+1},\cdots,$
then we have
$$
\text{$T_x$ is subnormal}\Longleftrightarrow\text{$T_x$ is 2-hyponormal}
\Longleftrightarrow
\cases
x=\alpha_j^2\ \ &\text{if}\ j\ge 1;\\
x\le a &\text{if}\ j=0.
\endcases
$$
\endproclaim

\demo{Proof}
Since $\alpha$ is recursively generated by $\sqrt{a},\sqrt{b},\sqrt{c}$,
we have that $\alpha_0^2=a,\ \alpha_1^2=b,\ \alpha_2^2=c$,
$$
\alpha_3^2=\frac{b(c^2-2ac+ab)}{c(b-a)},\quad\text{and}\quad
\alpha_4^2=\frac{bc^3-4abc^2+2ab^2c+a^2bc-a^2b^2+a^2c^2}{(b-a)(c^2-2ac+ab)}.
\tag 3.4
$$

{\it Case 1} ($j=0$): It is evident that $T_x$ is subnormal if and only if
$x\le a$.
For 2-hyponormality observe by (2.12) that $T_x$ is 2-hyponormal if
and only if
$$
\text{det}\pmatrix
1&x&bx\\ x&bx&bcx\\ bx&bcx&\alpha_3^2bcx\endpmatrix\ge 0,
$$
or equivalently, $x\le a$.

{\it Case 2} ($j\ge 1$): Without loss of generality we may assume
that $j=1$ and $a=1$. Thus $\alpha_1=\sqrt{x}$. Then by Theorem 2.1,
$T_x$ is subnormal if and only if $x=b$. On the other hand, by (2.12),
$T_x$ is 2-hyponormal if and only if
$$
\text{det}\pmatrix
1&1&x\\ 1&x&cx\\
x&cx&\alpha_3^2cx\endpmatrix\ge 0\quad\text{and}\quad \text{det}\pmatrix
1&x&cx\\ x&cx&\alpha_3^2cx\\
cx&\alpha_3^2cx&\alpha_3^2\alpha_4^2cx\endpmatrix\ge 0.
$$
Thus a direct calculation with the specific forms of $\alpha_3,\alpha_4$
given in (3.4)
shows that $T_x$ is 2-hyponormal if and only if
$(x-b)\left(x-\frac{b(c^2-2c+b)}{b-1}\right)\le 0$ and $x\le b$.
Since $b\le \frac{b(c^2-2c+b)}{b-1}$, it follows that $T_x$ is
2-hyponormal if and only if $x=b$. This completes the proof.
\hfill$\square$
\enddemo
\bigskip

\head 4. Proof of Theorem 2.2
\endhead
\bigskip

With the notation in  (2.6), we let
$$
p_n:=u_n\,v_{n+1}-w_n\qquad\text{($n\ge 0$)}.
$$
We then have:

\proclaim{Lemma 4.1} If $\alpha\equiv\{\alpha_n\}_{n=0}^\infty$ is
a strictly increasing weight sequence then the following
statements are equivalent:
\roster
\item"(i)" $W_\alpha$ is 2-hyponormal;
\item"(ii)"
$\alpha_{n+1}^2(u_{n+1}+u_{n+2})^2 \le
u_{n+1}v_{n+2}$
\quad($n\ge 0$);
\item"(iii)"
$\frac{\alpha_{n}^2}{\alpha_{n+2}^2}\frac{u_{n+2}}
{u_{n+3}}\ \le\ \frac{u_{n+1}}
{u_{n+2}}$\quad($n\ge 0$);
\item"(iv)" $p_n\ge 0$\quad($n\ge 0$).
\endroster
\endproclaim

\demo{Proof} This follows from a straightforward calculation.
\hfill$\square$
\enddemo
\bigskip

\demo{Proof of Theorem 2.2}
If $\alpha$ is not strictly increasing then $\alpha$ is flat,
by the argument of [7, Corollary 6], i.e.,
$\alpha_0=\alpha_1=\alpha_2=\cdots$. Then
$$
D_n(s)=\left(\smallmatrix \alpha_0^2+|s|^2\alpha_0^4& \bar s \alpha_0^3\\
s\alpha_0^3&|s|^2\alpha_0^4\endsmallmatrix\right)\oplus 0_\infty
\tag 4.1
$$
(cf. (2.5)), so that (2.9) is evident.
Thus we may assume that $\alpha$ is
strictly increasing, so that $u_n>0,\ v_n>0$ and $w_n>0$ for all
$n\ge 0$. Recall that if we write $d_n(t):=\sum_{i=0}^{n+1} c(n,i)t^i$
then the $c(n,i)$'s satisfy the following recursive formulas (cf. (2.8)):
$$
c(n+2,i)=u_{n+2}\,c(n+1,i)+v_{n+2}\,c(n+1,i-1)-w_{n+1}\,c(n,i-1)\ \
(n\ge 0,\ 1\le i\le n).
\tag 4.2
$$
Also, $c(n,n+1)=v_0\cdots v_n$ (again by (2.8)) and
$p_n:=u_n v_{n+1}-w_n\ge 0$ ($n\ge 0$), by Lemma 4.1.
A straightforward calculation shows that
$$
\align
d_0(t)&=u_0+v_0\,t;
\tag 4.3\\
d_1(t)&=u_0u_1+(v_0u_1+p_0)\,t+v_0v_1\,t^2;\\
d_2(t)&=u_0u_1u_2+(v_0u_1u_2+u_0p_1+u_2p_0)\,t
+(v_0v_1u_2+v_0p_1+v_2p_0)\,t^2+v_0v_1v_2\,t^3.
\endalign
$$
Evidently,
$$
c(n,\,i)\ge 0\qquad\text{($0\le n\le 2,\ 0\le i\le n+1$)}.
\tag 4.4
$$
Define
$$
\beta(n,i):=c(n,i)-v_0\cdots v_{i-1}u_i\cdots u_n\qquad
\text{($n\ge 1,\ 1\le i\le n$)}.
$$
For every $n\ge 1$, we now have
$$
c(n,i)=\cases
u_0\cdots u_n\ge 0 &(i=0)\\
v_0\cdots v_{i-1}u_i\cdots u_n\ +\ \beta(n,i)\quad&(1\le i\le n)\\
v_0\cdots v_n\ge 0 &(i=n+1).
\endcases
\tag 4.5
$$
For notational convenience we let $\beta(n,0):=0$ for every $n\ge 0$.
\medskip

\proclaim{Claim 1}
For $n\ge 1$,
$$
c(n,\,n)\ge u_n\,c(n-1,\,n)\ge 0.
\tag 4.6
$$
\endproclaim

\noindent
{\it Proof of Claim 1}. We use mathematical induction.
For $n=1$,
$$
c(1,1)=v_0u_1+p_0\ \ge\ u_1\,c(0,1)\ge 0,
$$
and
$$
\align
c(n+1,n+1)
&=u_{n+1}\,c(n,n+1)+v_{n+1}\,c(n,n)-w_{n}c(n-1,n)\\
&\ge u_{n+1}\,c(n,n+1)+v_{n+1}\,u_{n}c(n-1,n)-w_{n}\,c(n-1,n)
\quad\text{(by inductive hypothesis)}\\
&=u_{n+1}\,c(n,n+1)+p_{n}\,c(n-1,n)\\
&\ge u_{n+1}\,c(n,n+1),
\endalign
$$
which proves Claim 1.
\medskip

\proclaim{Claim 2} For $n\ge 2$,
$$
\beta(n,\,i)\ge u_{n}\,\beta(n-1,i)\ge 0\qquad (0\le i\le n-1).
\tag 4.7
$$
\endproclaim

\noindent
{\it Proof of Claim 2}.
We use mathematical induction. If $n=2$ and $i=0$,
this is trivial. Also,
$$
\beta(2,1)=u_0\,p_1+u_2\,p_0=u_0\,p_1+u_2\,\beta(1,1)\ge
u_2\,\beta(1,1)\ge 0.
$$
Assume that (4.7) holds. We shall prove that
$$
\beta(n+1,\,i)\ge u_{n+1}\,\beta(n,i)\ge 0\qquad (0\le i\le n).
$$
For,
$$
\align
\beta(n+1&,i)+v_0\cdots v_{i-1}u_i\cdots u_{n+1}=c(n+1,i)\qquad
\text{(by (4.2))}\\
&=u_{n+1} c(n,i)+v_{n+1} c(n,i-1)-w_n c(n-1,i-1)\\
&=u_{n+1}\biggl(\beta(n,i)+v_0\cdots v_{i-1}u_i\cdots u_n\biggr)\\
&\qquad +v_{n+1}\biggl(\beta(n,i-1)+v_0\cdots v_{i-2}u_{i-1}\cdots
u_{n}\biggr)\\
&\qquad\qquad -w_n\biggl(\beta(n-1,i-1)+v_0\cdots v_{i-2}
u_{i-1}\cdots u_{n-1}\biggr),
\endalign
$$
so that
$$
\align
\beta(n+1,i)
&= u_{n+1}\beta(n,i)+v_{n+1}\beta(n,i-1)-w_n\beta(n-1,i-1)\\
&\quad +
     v_0\cdots v_{i-2}u_{i-1}\cdots u_{n-1}\, (u_n v_{n+1}-w_n)\\
&= u_{n+1}\beta(n,i)+v_{n+1}\beta(n,i-1)-w_n\beta(n-1,i-1)+
     (v_0\cdots v_{i-2}u_{i-1}\cdots u_{n-1})\,p_n\\
&\ge u_{n+1}\beta(n,i)+v_{n+1}u_n\beta(n-1,i-1)-w_n\beta(n-1,i-1)\\
&\qquad\text{(by the
inductive hypothesis and Lemma 4.1;}\\
&\qquad\ \text{observe that $i-1\le n-1$,
so (4.7) applies)}\\
&=u_{n+1}\beta(n,i)+p_n\,\beta(n-1,i-1)\\
&\ge u_{n+1}\,\beta(n,i),
\endalign
$$
which proves Claim 2.

By Claim 2 and (4.5), we can see that $c(n,i)\ge 0$ for all $n\ge 0$ and
$1\le i\le n-1$. Therefore (4.4), (4.5), Claim 1 and Claim 2 imply
$$
c(n,i)\ \ge\ v_0\cdots v_{i-1}u_i\cdots u_n\qquad\text{($n\ge 0,\
0\le i\le n+1$)}.
$$
This completes the proof.
\hfill$\square$
\enddemo
\bigskip

\head 5. Proof of Theorem 2.3
\endhead
\bigskip

To prove Theorem 2.3 we need:

\proclaim{Lemma 5.1 ([15, Lemma 2.3])}
Let $\alpha\equiv\{\alpha_n\}_{n=0}^\infty$ be a strictly
increasing weight sequence. If $W_\alpha$ is 2-hyponormal
then the sequence of quotients
$$
\Theta_n:=\frac{u_{n+1}}{u_{n+2}}\qquad
(n\ge 0)
\tag 5.1
$$
is bounded away from $0$ and from $\infty$. More precisely,
$$
1\le \Theta_n\le \frac{u_1}{u_2}
\left(\frac{||W_\alpha||^2}{\alpha_0\alpha_1}\right)^2\quad\text{
for sufficiently large $n$}.
\tag 5.2
$$
In particular, $\{u_n\}_{n=0}^\infty$ is eventually decreasing.
\endproclaim
\medskip

\demo{Proof of Theorem 2.3}
By Theorem 2.2, $W_\alpha$ is {\it strictly} positively quadratically
hyponormal, in the sense that
all coefficients of $d_n(t)$ are {\it positive} for all $n\ge 0$. Note that
finite rank perturbations of $\alpha$ affect
a finite number of values of $u_n,\ v_n$ and $w_n$.
More concretely, if $\alpha^\prime$ is a perturbation of $\alpha$
in the weights $\{\alpha_0,\cdots, \alpha_N\}$, then
$u_n,\ v_n$, $w_n$ and $p_n$ are invariant under $\alpha^\prime$ for
$n\ge N+3$.
In particular, $p_n\ge 0$ for $n\ge N+3$.
\medskip

\proclaim{Claim 1} For $n\ge 3,\ 0\le i\le n+1$,
$$
\align
c(n,i)=&u_n\,c(n-1,i)+p_{n-1}\,c(n-2,i-1)
+\sum_{k=4}^n p_{k-2}\left(\prod_{j=k}^n v_j\right)c(k-3,\,i-n+k-2)\\
&+v_n\cdots v_3\,\rho_{i-n+1},
\tag 5.3
\endalign
$$
where
$$
\rho_{i-n+1}=\cases
0&(i<n-1)\\
u_0p_1&(i=n-1)\\
v_0p_1+v_2p_0\quad&(i=n)\\
v_0v_1v_2&(i=n+1)
\endcases
$$
\endproclaim\noindent
(cf. [12, Proof of Theorem 4.3]).
\medskip

\noindent
{\it Proof of Claim 1.} We use induction.
For $n=3,\ 0\le i\le 4$,
$$
\align
c(3,i)
&=u_3\,c(2,\,i)+v_3\,c(2,\,i-1)-w_2\,c(1,\,i-1)\\
&=u_3\,c(2,i)+v_3\biggl(u_2\,c(1,i-1)+v_2\,c(1,\,i-2)-w_1\,c(0,\,i-2)\biggr)
    -w_2\,c(1,i-1)\\
&=u_3\,c(2,i)+p_2\,c(1,i-1)+v_3\biggl(v_2\,c(1,i-2)-w_1\,c(0,i-2)\biggr)\\
&=u_3\,c(2,i)+p_2\,c(1,i-1)+v_3\,\rho_{i-2},
\endalign
$$
where by (4.3),
$$
\rho_{i-2}=\cases
0&(i<2)\\
u_0p_1&(i=2)\\
v_0p_1+v_2p_0\quad&(i=3)\\
v_0v_1v_2&(i=4).
\endcases
$$
Now,
$$
\align
c(n+1,i)
&=u_{n+1}c(n,i)+v_{n+1}c(n,i-1)-w_n c(n-1,i-1)\\
&=u_{n+1}c(n,i)+v_{n+1}\biggl(u_n c(n-1,i-1)+p_{n-1}c(n-2,i-2)\\
&\ \ +
\sum_{k=4}^n p_{k-2}\left(\prod_{j=k}^n v_j\right)c(k-3,i-n+k-3)+
v_n\cdots v_3 \rho_{i-n}\biggr)-w_n\,c(n-1,i-1)\\
&=u_{n+1}c(n,i)+p_n c(n-1,i-1)+v_{n+1}p_{n-1}c(n-2,i-2)\\
&\ \ +
v_{n+1}\sum_{k=4}^{n}p_{k-2}\left(\prod_{j=k}^{n} v_j\right) c(k-3,i-n+k-3)
+v_{n+1}\cdots v_3\rho_{i-n}\\
&\qquad \text{(by inductive hypothesis)}\\
&=u_{n+1}c(n,i)+p_n c(n-1,i-1)+
\sum_{k=4}^{n+1}p_{k-2}\left(\prod_{j=k}^{n+1} v_j\right) c(k-3,i-n+k-3)\\
&\ \ +v_{n+1}\cdots v_3\rho_{i-n},
\endalign
$$
which proves Claim 1.
\medskip

Write $u_n^\prime,\ v_n^\prime, w_n^\prime, p_n^\prime, \rho_n^\prime,$
and $c^\prime(\cdot,\cdot)$ for the entities corresponding to $\alpha^\prime$.
If $p_n>0$ for every $n=0,\cdots, N+2$, then in view of Claim 1,
we can choose a small perturbation such that $p_n^\prime>0$ ($0\le n\le N+2$)
and therefore
$c^\prime(n,i)>0$ for all $n\ge 0$ and $0\le i\le n+1$, which
implies that $W_{\alpha^\prime}$ is also positively quadratically hyponormal.
If instead $p_n=0$ for some $n=0,\cdots, N+2$,
careful inspection of (5.3) reveals that
without loss of generality we may assume $p_0=\cdots=p_{N+2}=0$.
By Theorem 2.2, we have that for a sufficiently small perturbation
$\alpha^\prime$ of $\alpha$,
$$
c^\prime(n,i)>0\ \ (0\le n\le N+2,\ 0\le i\le n+1)
\quad\text{and}\quad c^\prime(n,n+1)>0\ \ (n\ge 0).
\tag 5.4
$$
Write
$$
k_n:=\frac{v_n}{u_n}\qquad(n=2,3,\cdots).
$$
\medskip

\proclaim{Claim 2} $\{k_n\}_{n=2}^\infty$ is bounded.
\endproclaim

\noindent
{\it Proof of Claim 2}.
Observe that
$$
\align
k_n=\frac{v_n}{u_n}
&=\frac{\alpha_n^2\alpha_{n+1}^2-\alpha_{n-1}^2\alpha_{n-2}^2}
{\alpha_n^2-\alpha_{n-1}^2}\\
&=\alpha_n^2+\alpha_{n-1}^2+\alpha_n^2
\frac{\alpha_{n+1}^2-\alpha_n^2}{\alpha_n^2-\alpha_{n-1}^2}+\alpha_{n-1}^2
\frac{\alpha_{n-1}^2-\alpha_{n-2}^2}{\alpha_n^2-\alpha_{n-1}^2}.
\tag 5.5
\endalign
$$
Therefore if $W_\alpha$ is 2-hyponormal then by Lemma 5.1,
the sequences
$$
\left\{\frac{\alpha_{n+1}^2-\alpha_n^2}{\alpha_n^2-\alpha_{n-1}^2}
\right\}_{n=2}^\infty
\quad\text{and}\quad
\left\{\frac{\alpha_{n-1}^2-\alpha_{n-2}^2}{\alpha_n^2-\alpha_{n-1}^2}
\right\}_{n=2}^\infty
$$
are both bounded, so that
$\{k_n\}_{n=2}^\infty$ is bounded. This proves Claim 2.
\medskip

Write  $k:=\sup_n k_n$.
Without loss of generality
we assume $k<1$ (this is possible from the observation that
$c\alpha$ induces $\{c^2 k_n\}$).
Choose a sufficiently small perturbation $\alpha^\prime$
of $\alpha$ such that if we let
$$
h:=\sup_{\Sb 0\le \ell\le N+2\\
0\le m\le 1\endSb}
\left|\sum_{k=4}^{N+4}p_{k-2}^\prime\left(\prod_{j=k}^{N+3}v_j^\prime\right)
c^\prime(k-3,\,\ell)+v_{N+3}^\prime\cdots v_3^\prime\,\rho_m^\prime
\right|
\tag 5.6
$$
then
$$
c^\prime(N+3,\,i)-\frac{1}{1-k}h\ >\ 0\qquad\text{($0\le i\le N+3$)}
\tag 5.7
$$
(this is always possible because by Theorem 2.2, we can choose a
sufficiently small $|p_i^\prime|$ such that
$$
c^\prime(N+3,i)>v_0\cdots v_{i-1}u_i\cdots u_{N+3}-\epsilon
\quad\text{and}\quad
|h|<(1-k)\bigl(v_0\cdots v_{i-1}u_i\cdots u_{N+3}-\epsilon\bigr)
$$
for any small $\epsilon>0$).
\medskip

\proclaim{Claim 3}
For $j\ge 4$ and $0\le i\le N+j$,
$$
c^\prime(N+j,\,i)\ \ge\ u_{N+j}\cdots u_{N+4}
\left(c^\prime(N+3,\ i)-\sum_{n=1}^{j-3} k^n\,h\right).
\tag 5.8
$$
\endproclaim

\noindent
{\it Proof of Claim 3}.
We use induction. If $j=4$
then by Claim 1 and (5.6),
$$
\align
c^\prime(N+4,\,i)
&= u_{N+4}^\prime c^\prime(N+3,i)+p_{N+3}^\prime c^\prime
(N+2,i-1)\\
&\quad
+v_{N+4}^\prime\sum_{k=4}^{N+4} p_{k-2}^\prime
\left(\prod_{j=k}^{N+3} v_j^\prime\right)c^\prime(k-3,i-N+k-6)
+v_{N+4}^\prime\cdots v_3^\prime\rho^\prime_{i-(N+3)}\\
&\ge u_{N+4}^\prime c^\prime(N+3,i)+p_{N+3}^\prime c^\prime
(N+2,i-1)-v_{N+4}^\prime h\\
&\ge u_{N+4}\bigl(c^\prime(N+3,i)-k_{N+4}h\bigr)\\
&\ge u_{N+4}\bigl(c^\prime(N+3,i)-k\,h\bigr),
\endalign
$$
because $u_{N+4}^\prime=u_{N+4},\ v_{N+4}^\prime=v_{N+4}$
and $p_{N+3}^\prime=p_{N+3}\ge 0$.
Now suppose (5.8) holds for some $j\ge 4$.
By Claim 1, we have that for $j\ge 4$,
$$
\align
c^\prime (N+j+1,i)
&=u_{N+j+1}^\prime c^\prime (N+j,\,i)+p_{N+j}^\prime
c(N+j-1,\,i-1)\\
&\quad
+\sum_{k=4}^{N+j+1} p_{k-2}^\prime \left(
\prod_{j=k}^{N+j+1} v_j^\prime\right)c^\prime(k-3, i-N+k-j-3)
+v_{N+j+1}^\prime\cdots v_3^\prime \rho_{i-(N+j)}^\prime\\
&=u_{N+j+1}^\prime c^\prime (N+j,\,i)+p_{N+j}^\prime
c(N+j-1,\,i-1)\\
&\quad
+\sum_{k=N+5}^{N+j+1} p_{k-2}^\prime \left(
\prod_{j=k}^{N+j+1} v_j^\prime\right)c^\prime(k-3, i-N+k-j-3)\\
&\quad
+\sum_{k=4}^{N+4} p_{k-2}^\prime \left(
\prod_{j=k}^{N+j+1} v_j^\prime\right)c^\prime(k-3, i-N+k-j-3)
+v_{N+j+1}^\prime\cdots v_3^\prime \rho_{i-(N+j)}^\prime.
\endalign
$$
Since $p_n^\prime=p_n>0$ for $n\ge N+3$ and
$c^\prime(n,\ell)>0$ for $0\le n\le N+j$ by the inductive hypothesis,
it follows that
$$
p_{N+j}^\prime
c(N+j-1,\,i-1)+\sum_{k=N+5}^{N+j+1} p_{k-2}^\prime \left(
\prod_{j=k}^{N+j+1} v_j^\prime\right)c^\prime(k-3, i-N+k-j-3)
\ge 0.
\tag 5.9
$$
By inductive hypothesis and (5.9),
$$
\align
&c^\prime (N+j+1,i)\\
&\ge
u_{N+j+1}^\prime c^\prime (N+j,\,i)
+\sum_{k=4}^{N+4} p_{k-2}^\prime \left(
\prod_{j=k}^{N+j+1} v_j^\prime\right)c^\prime(k-3, i-N+k-j-3)
+v_{N+j+1}^\prime\cdots v_3^\prime \rho_{i-(N+j)}^\prime\\
&\ge
u_{N+j+1} u_{N+j}\cdots u_{N+4}
\left(c^\prime(N+3,i)-\sum_{n=1}^{j-3} k^n h\right)\\
&\ \
+v_{N+j+1}v_{N+j}\cdots v_{N+4}
\left(
\sum_{k=4}^{N+4} p_{k-2}^\prime \left(
\prod_{j=k}^{N+3} v_j^\prime\right)c^\prime(k-3, i-N+k-j-3)
+v_{N+3}^\prime\cdots v_3^\prime \rho_{i-(N+j)}^\prime
\right)\\
&\ge
u_{N+j+1} u_{N+j}\cdots u_{N+4}
\left(c^\prime(N+3,i)-\sum_{n=1}^{j-3} k^n h\right)
-v_{N+j+1}v_{N+j}\cdots v_{N+4}\,h\\
&=
u_{N+j+1} u_{N+j}\cdots u_{N+4}
\left(c^\prime(N+3,i)-\sum_{n=1}^{j-3} k^n h
-k_{N+j+1}k_{N+j}\cdots k_{N+4}\,h\right)\\
&\ge
u_{N+j+1} u_{N+j}\cdots u_{N+4}
\left(c^\prime(N+3,i)-\sum_{n=1}^{j-2} k^n h\right),
\endalign
$$
which proves Claim 3.
\medskip

Since $\sum_{n=1}^j k^n<\frac{1}{1-k}$ for every $j>1$, it follows from
Claim 3 and (5.7) that
$$
c^\prime(N+j,\,i)>0\quad\text{for}\ j\ge 4\ \text{and}\ 0\le i\le N+j.
\tag 5.10
$$
It thus follows from (5.4) and (5.10) that
$c^\prime(n,i)>0$ for every $n\ge 0$ and $0\le i\le n+1$.
Therefore
$W_{\alpha^\prime}$ is also positively quadratically hyponormal.
This completes the proof.
\hfill$\square$
\enddemo
\bigskip

\proclaim{Corollary 5.2} Let $W_\alpha$ be a weighted shift such that
$\alpha_{j-1}<\alpha_j$ for some $j\ge 1$, and let
$T_x$ be the weighted shift
with  weight sequence
$$
\alpha_x:\alpha_0,\cdots,\alpha_{j-1},x,\alpha_{j+1},\cdots.
$$
Then $\{x: T_x\ \text{is 2-hyponormal}\}$ is a proper closed
subset of $\{x:\ T_x\ \text{is quadratically hyponormal}\}$
whenever the
latter set is non-empty.
\endproclaim

\demo{Proof}
Write
$$
H_2:=\{x: T_x\ \text{is 2-hyponormal}\}.
$$
Without loss of generality,
we can assume that
$H_2$ is non-empty, and that $j=1$.
Recall that a 2-hyponormal weighted shift with two equal weights
is of the form
$\alpha_0=\alpha_1=\alpha_2=\cdots$ or $\alpha_0<\alpha_1=\alpha_2
=\alpha_3=\cdots$.
Let $x_m:=\inf\,H_2$. By Proposition 6.7 below,
$T_{x_m}$ is hyponormal. Then $x_m>\alpha_0$. By assumption,
$x_m<\alpha_2$. Thus $\alpha_0,x_m,\alpha_2,\alpha_3,\cdots$ is
strictly increasing.
Now we apply Theorem 2.3 to obtain
$x^\prime$ such that $\alpha_0<x^\prime<x_m$ and $T_{x^\prime}$
is quadratically hyponormal. However
$T_{x^\prime}$ is not 2-hyponormal by the definition of $x_m$.
The proof is complete.
\hfill$\square$
\enddemo
\bigskip

The following question arises naturally:

\proclaim{Question 5.3}
Let $\alpha$ be a strictly increasing weight sequence
and let $k\ge 3$.
If $W_\alpha$ is a $k$-hyponormal weighted shift,
does it follow that
$W_\alpha$ is weakly $k$-hyponormal under a small perturbation of
the weight sequence\,?
\endproclaim
\bigskip

\head 6. Other Related Results
\endhead
\bigskip

\noindent
{\bf \S 6.1 Subnormal Extensions}
\medskip

\noindent
Let $\alpha:\alpha_0,\alpha_1,\cdots$
be a weight sequence, let $x_i>0$ for $1\le i\le n$, and let
$(x_n,\cdots x_1)\alpha:x_n,\cdots,x_1,\alpha_0,\alpha_1,\cdots$
be the augmented weight sequence. We say that $W_{(x_n,\cdots,x_1)\alpha}$
is an {\it extension} (or {\it $n$-step extension}) of $W_\alpha$.
Observe that
$$
W_{(x_n,\cdots,x_1)\alpha}\vert_{\bigvee\{e_n,e_{n+1},\cdots\}}
\cong W_\alpha.
$$
The hypothesis $F\ne c\,P_{\{e_0\}}$ in Theorem 2.1 is essential.
Indeed, there exist infinitely many one-step subnormal extension
of a subnormal
weighted shift whenever one such extension exists. Recall
([7, Proposition 8]) that if $W_\alpha$ is a weighted shift
whose restriction to $\bigvee\{e_1,e_2,\cdots\}$ is subnormal with
associated measure $\mu$, then $W_\alpha$ is subnormal if and only if
\roster
\item"(i)" $\frac{1}{t}\in L^1(\mu)$;
\item"(ii)" $\alpha_0^2\le \left(||\frac{1}{t}||_{L^1(\mu)}\right)^{-1}$.
\endroster
Also note that there may not exist any one-step subnormal extension of
the subnormal weighted shift: for example, if $W_\alpha$ is the Bergman
shift then the corresponding Berger measure is $\mu(t)=t$, and hence
$\frac{1}{t}$ is not integrable with respect to $\mu$; therefore $W_\alpha$
does not admit any subnormal extension.
A similar situation arises when $\mu$ has an atom at $\{0\}$.
\smallskip

More generally we have:

\proclaim{Theorem 6.1 (Subnormal Extensions)}
Let $W_\alpha$ be a subnormal weighted shift with weights
$\alpha:\alpha_0,\alpha_1,\cdots$ and let $\mu$ be the corresponding
Berger measure. Then $W_{(x_n,\cdots,x_1)\alpha}$ is subnormal if and
only if
\roster
\item"(i)" $\frac{1}{t^n}\in L^1(\mu)$;
\item"(ii)" $x_j=\left(\frac{||\frac{1}{t^{j-1}}||_{L^1(\mu)}}{||
\frac{1}{t^j}||_{L^1(\mu)}}\right)^{\frac{1}{2}}$\quad for $1\le j
\le n-1$;
\item"(iii)" $x_n\le \left(\frac{||\frac{1}{t^{n-1}}||_{L^1(\mu)}}
{||\frac{1}{t^n}||_{L^1(\mu)}}\right)^{\frac{1}{2}}$.
\endroster
In particular, if we put
$$
S:=\{(x_1,\cdots,x_n)\in \Bbb{R}^n: W_{(x_n,\cdots,x_1)\alpha}\ \
\text{is subnormal}\}
$$
then either $S=\emptyset$ or $S$ is a line segment in $\Bbb{R}^n$.
\endproclaim

\demo{Proof}
Write $W_j:=W_{(x_n,\cdots,x_1)\alpha}\vert_{\bigvee\{e_{n-j},e_{n-j+1},
\cdots\}}$ ($1\le j\le n$) and hence $W_n=W_{(x_n,\cdots,x_1)\alpha}$.
By the argument used to establish  (3.2) we have that
$W_1$ is subnormal with associated measure $\nu_1$ if and only if
\roster
\item"(i)" $\frac{1}{t}\in L^1(\mu)$;
\item"(ii)" $d\nu_1=\frac{x_1^2}{t}d\mu$, or equivalently, $x_1^2=
\left(\int_0^{||W_\alpha||^2}\frac{1}{t}\,d\mu(t)\right)^{-1}$.
\endroster
Inductively $W_{n-1}$ is subnormal with associated measure
$\nu_{n-1}$ if and only if
\roster
\item"(i)" $W_{n-2}$ is subnormal;
\item"(ii)" $\frac{1}{t^{n-1}}\in L^1(\mu)$;
\item"(iii)" $d\nu_{n-1}=\frac{x_{n-1}^2}{t}d\nu_{n-2}=\cdots
=\frac{x_{n-1}^2\cdots x_1^2}{t^{n-1}}d\mu$, or equivalently,
$x_{n-1}^2=\frac{\int_0^{||W_\alpha||^2}\frac{1}{t^{n-2}}\,d\mu(t)}
{\int_0^{||W_\alpha||^2}\frac{1}{t^{n-1}}\,d\mu(t)}$.
\endroster
Therefore $W_n$ is subnormal if and only if
\roster
\item"(i)" $W_{n-1}$ is subnormal;
\item"(ii)" $\frac{1}{t^n}\in L^1(\mu)$;
\item"(iii)" $x_n^2\le \left(\int_0^{||W_\alpha||^2} \frac{1}{t}
d\nu_{n-1}\right)^{-1}=
\left(\int_0^{||W_\alpha||^2} \frac{x_{n-1}^2\cdots x_1^2}{t^n}
d\mu(t)\right)^{-1}
=\frac{\int_0^{||W_\alpha||^2}\frac{1}{t^{n-1}}\,d\mu(t)}
{\int_0^{||W_\alpha||^2}\frac{1}{t^n}\,d\mu(t)}$.
\endroster
\hfill$\square$
\enddemo
\bigskip

\proclaim{Corollary 6.2}
If $W_\alpha$ is a subnormal weighted shift with associated measure
$\mu$, there exists an $n$-step subnormal extension of $W_\alpha$ if
and only if $\frac{1}{t^n}\in L^1(\mu)$.
\endproclaim
\bigskip

For the next result we refer to the notation in
(2.1) and (2.2).

\proclaim{Corollary 6.3}
A recursively generated subnormal shift with $\varphi_0\ne 0$
admits an $n$-step subnormal
extension for every $n\ge 1$.
\endproclaim

\demo{Proof}
The assumption about $\varphi_0$ implies that
the zeros of $g(t)$ are positive,
so that $s_0>0$. Thus for every $n\ge 1$, $\frac{1}{t^n}$
is integrable with respect to the corresponding Berger measure
$\mu=\rho_0\delta_{s_0}+\cdots+\rho_{r-1}\delta_{s_{r-1}}$.
By Corollary 6.2,
there exists an $n$-step subnormal extension.
\hfill$\square$
\enddemo
\bigskip

We need not expect that for arbitrary recursively generated
shifts,
2-hyponormality and
subnormality coincide as in Theorem 3.2.
For example, if $\alpha:\sqrt{\frac{1}{2}},\sqrt{x},(\sqrt{3},
\sqrt{\frac{10}{3}},\sqrt{\frac{17}{5}})^\wedge$ then by (2.12)
and Theorem 6.1,
\roster
\item"(i)" $T_x$ is 2-hyponormal $\Longleftrightarrow$ $4-\sqrt{6}\le x\le 2$;
\item"(ii)" $T_x$ is subnormal $\Longleftrightarrow$ $x=2$.
\endroster
A straightforward calculation shows, however, that
$T_x$ is {\it $3$-hyponormal} if and only if $x=2$; for,
$$
A(0;3):=\pmatrix
1&\frac{1}{2}&\frac{1}{2}x&\frac{3}{2}x\\
\frac{1}{2}&\frac{1}{2}x&\frac{3}{2}x&5x\\
\frac{1}{2}x&\frac{3}{2}x&5x&17x\\
\frac{3}{2}x&5x&17x&58x\endpmatrix\ \ge\ 0\ \Longleftrightarrow\ x=2.
$$
This behavior is typical of general recursively generated
weighted shifts: we show in [13] that
subnormality is equivalent to
$k$-hyponormality for some $k\ge 2$.
\bigskip

\noindent
{\bf \S 6-2 Convexity and Closedness}
\medskip

\noindent
Next, we will show that canonical
rank-one perturbations of $k$-hyponormal
weighted shifts which preserve
$k$-hyponormality form a convex set.
To see this we need an auxiliary result.

\proclaim{Lemma 6.4}
Let $I=\{1,\cdots,n\}\times \{1,\cdots,n\}$ and let $J$ be a symmetric
subset of $I$.
Let $A=(a_{ij})\in M_n(\Bbb{C})$ and let $C=(c_{ij})\in
M_n(\Bbb{C})$ be given by
$$
c_{ij}=\cases
c\,a_{ij}\quad &\text{if}\ \ (i,j)\in J\\
a_{ij}&\text{if}\ \ (i,j)\in I\setminus J
\endcases
\qquad (c>0).
$$
If $A$ and $C$ are positive semidefinite then $B=(b_{ij})\in M_n(\Bbb{C})$
defined by
$$
b_{ij}=\cases
b\,a_{ij}\quad &\text{if}\ \ (i,j)\in J\\
a_{ij}&\text{if}\ \ (i,j)\in I\setminus J\endcases
\qquad (b\in [1,c]\ \text{or}\ [c,1])
$$
is also positive semidefinite.
\endproclaim

\demo{Proof}
Without loss of generality we may assume
$c>1$. If $b=1$ or $b=c$ the assertion is trivial.
Thus we assume $1<b<c$.
The result is now a consequence of the following observation.
If $[D]_{(i,j)}$ denotes the ($i,j$)-entry of the matrix
$D$ then
$$
\align
\left[\frac{c-b}{c-1}\left(A+\frac{b-1}{c-b}C\right)\right]_{(i,j)}
&=\cases
\frac{c-b}{c-1}\left(1+\frac{b-1}{c-b}c\right)a_{ij}\quad
&\text{if}\ \  (i,j)\in J\\
\frac{c-b}{c-1}\left(1+\frac{b-1}{c-b}\right)a_{ij}
&\text{if}\ \ (i,j)\in I\setminus J\endcases\\
&=\cases b\,a_{ij}\quad&\text{if}\ \ (i,j)\in J\\
a_{ij}&\text{if}\ \ (i,j)\in I\setminus J\endcases\\
&=[B]_{(i,j)},
\endalign
$$
which is positive semidefinite because positive
semidefinite matrices in $M_n(\Bbb{C})$ form a cone.
\hfill$\square$
\enddemo
\bigskip

An immediate consequence of Lemma 6.4 is that positivity of a
matrix forms a convex set with respect to a fixed diagonal
location; i.e., if
$$
A_x=\pmatrix *&*&*\\*&x&*\\ *&*&*\endpmatrix
$$
then $\{x:A_x\ \text{is positive semidefinite}\}$ is convex.
\bigskip

We now have:

\proclaim{Theorem 6.5}
Let $\alpha=\{\alpha_n\}_{n=0}^\infty$ be a weight sequence,
let $k\ge 1$, and let $j\ge 0$. Define $\alpha^{(j)}(x):
\alpha_0,\cdots,\alpha_{j-1},x,\alpha_{j+1},\cdots$.
Assume $W_\alpha$ is $k$-hyponormal
and define
$$
\Omega_\alpha^{k,j}:=
\{ x:\ W_{\alpha^{(j)}(x)}\ \text{is $k$-hyponormal}\}.
$$
Then $\Omega_\alpha^{k,j}$ is a closed interval.
\endproclaim

\demo{Proof}
Suppose $x_1,x_2\in \Omega_\alpha^{k,j}$ with $x_1<x_2$.
Then by (2.11), the $(k+1)\times (k+1)$ Hankel matrix
$$
A_{x_i}(n;k):=\pmatrix
\gamma_n&\gamma_{n+1}&\hdots&\gamma_{n+k}\\
\gamma_{n+1}&\gamma_{n+2}&\hdots&\gamma_{n+k+1}\\
\vdots&\vdots&&\vdots\\
\gamma_{n+k}&\gamma_{n+k+1}&\hdots&\gamma_{n+2k}\endpmatrix
\qquad(n\ge 0;\ i=1,2)
$$
is positive, where $A_{x_i}$ corresponds to $\alpha^{(j)}(x_i)$.
We must show that $t x_1+(1-t)x_2\in \Omega_\alpha^{k,j}$
($0<t<1$), i.e.,
$$
A_{t x_1+(1-t)x_2}(n;k)\ \ge\ 0\quad(n\ge 0,\ 0<t<1).
$$
Observe that it suffices to establish the positivity of the $2k$
Hankel matrices corresponding to
$\alpha^{(j)}(t x_1+(1-t)x_2)$
such that $t x_1+(1-t) x_2$ appears as a factor in at least one
entry but not in every entry.
A moment's thought reveals that without loss of generality
we may assume $j=2k$.
Observe that
$$
A_{z_1}(n;k)-A_{z_2}(n;k)=(z_1^2-z_2^2)\,H(n;k)
$$
for some Hankel matrix $H(n;k)$.
For notational convenience, we abbreviate $A_z(n;k)$ as $A_z$.
Then
$$
A_{t x_1+(1-t) x_2}=
\cases
t^2 A_{x_1}+(1-t)^2 A_{x_2}+2t(1-t)A_{\sqrt{x_1x_2}}
\quad&\text{for}\ 0\le n\le 2k\\
\left(t+(1-t)\frac{x_2}{x_1}\right)^2 A_{x_1}&\text{for}\ n\ge 2k+1.
\endcases
$$
Since $A_{x_1}\ge 0,\ A_{x_2}\ge 0$ and $A_{\sqrt{x_1x_2}}$ have the
form described by Lemma 6.4 and since $x_1<\sqrt{x_1x_2}<x_2$ it follows
from Lemma 6.4 that $A_{\sqrt{x_1x_2}}\ge 0$. Thus evidently,
$A_{t x_1+(1-t)x_2}\ge 0$, and therefore
$t x_1+(1-t)x_2\in \Omega_\alpha^{k,j}$.
This shows that $\Omega_\alpha^{k,j}$ is an interval.
The closedness of the interval follows from
Proposition 6.7 below.
\hfill$\square$
\enddemo
\bigskip

In [17] and [18], it was shown that there exists a non-subnormal
polynomially hyponormal  operator.
Also in [22], it was shown that there exists a non-subnormal
polynomially hyponormal operator
if and only if there exists one which is also a
weighted shift. However, no concrete weighted
shift has yet been found. As a strategy for finding
such a shift, we would like to suggest the following:

\proclaim{Question 6.6}
Does it follow that the polynomial hyponormality of the weighted shift is stable
under small perturbations of the weight sequence\,?
\endproclaim
\medskip

If the answer to Question 6.6 were affirmative then we would easily
find a polynomially hyponormal non-subnormal (even non-2-hyponormal)
weighted shift; for example, if
$$
\alpha: 1,\sqrt{x},(\sqrt{3},\sqrt{\frac{10}{3}},\sqrt{\frac{17}{5}})^\wedge
$$
and $T_x$ is the weighted shift associated with $\alpha$, then by Theorem 3.2,
$T_x$ is subnormal
$\Leftrightarrow$ $x=2$, whereas $T_x$ is polynomially hyponormal
$\Leftrightarrow$ $2-\delta_1<x<2+\delta_2$ for some $\delta_1,\delta_2>0$
provided the answer to Question 6.6 is yes; therefore for sufficiently
small $\epsilon>0$,
$$
\alpha_\epsilon: 1, \sqrt{2+\epsilon},(\sqrt{3},
\sqrt{\frac{10}{3}},\sqrt{\frac{17}{5}})^\wedge
$$
{\it would} induce a non-2-hyponormal polynomially hyponormal
weighted shift.
\bigskip

The answer to Question 6.6 for weak $k$-hyponormality is negative.
In fact we have:

\proclaim{Proposition 6.7}
\roster
\item"(i)" The set of $k$-hyponormal operators is $sot$-closed.
\item"(ii)" The set of weakly $k$-hyponormal operators is $sot$-closed.
\endroster
\endproclaim

\demo{Proof}
Suppose $T_\eta\in \Cal{L(H)}$ and $T_\eta\rightarrow T$ in $sot$. Then,
by the Uniform Boundedness Principle, $\{||T_\eta||\}_\eta$ is bounded. Thus
$T_\eta^{*i}T_\eta^j\rightarrow T^{*i}T^j$ in $sot$ for every $i,j$,
so that $M_k(T_\eta)\rightarrow M_k(T)$ in $sot$
(where $M_k(T)$ is as in (1.2)).
(i) In this case $M_k(T_\eta)\ge 0$ for all $\eta$, so $M_k(T)\ge 0$,
i.e., $T$ is $k$-hyponormal.

(ii) Here, $M_k(T_\eta)$ is weakly positive for all $\eta$.
By (1.3), $M_k(T)$ is also weakly positive, i.e., $T$ is weakly $k$-hyponormal.
\hfill$\square$
\enddemo

\vskip 2pc

\enddocument
\end